\documentclass[11pt,reqno]{amsart}
\usepackage{enumerate, latexsym, amsmath, amsfonts, amssymb, amsthm, color}
\def\pmod #1{\ ({\rm{mod}}\ #1)}
\def\Z{\Bbb Z}

\def\Q{\Bbb Q}

\def\R{\Bbb R}

\def\l{\left}
\def\r{\right}
\def\bg{\bigg}
\def\({\bg(}
\def\){\bg)}
\def\t{\text}
\def\f{\frac}

\def\per{{\rm per}}

\def\ls{\leqslant}

\def\se {\subseteq}
\def\sm{\setminus}
\def\bi{\binom}

\def\ve{\varepsilon}

\def\eq{\equiv}

\def\da{\delta}

\def\Proof{\noindent{\it Proof}}

\theoremstyle{plain}
\newtheorem{theorem}{Theorem}

\newtheorem{lemma}{Lemma}
\newtheorem{corollary}{Corollary}
\newtheorem{conjecture}{Conjecture}
\theoremstyle{definition}

\theoremstyle{remark}
\newtheorem{remark}{Remark}

\allowdisplaybreaks

 \vspace{4mm}

\begin{document}

\hbox{Preprint}
\medskip

\title
[{Some determinants involving quadratic residues modulo primes}]
{Some determinants involving \\ quadratic residues modulo primes}

\author
[Zhi-Wei Sun] {Zhi-Wei Sun}

\address{School of Mathematics, Nanjing
University, Nanjing 210093, People's Republic of China}
\email{zwsun@nju.edu.cn}

\keywords{Determinants, Legendre symbols, quadratic residues modulo primes, the tangent function.
\newline \indent Supported
by the National Natural Science Foundation of China (grant 12371004).}

\begin{abstract}
In this paper we evaluate several determinants involving quadratic residues modulo primes.
For example, for any prime $p>3$ with $p\equiv3\pmod4$ and $a,b\in\mathbb Z$ with $p\nmid ab$, we prove that
 $$\det\left[1+\tan\pi\frac{aj^2+bk^2}p\right]_{1\ls j,k\ls\f{p-1}2}=\begin{cases}-2^{(p-1)/2}p^{(p-3)/4}&\text{if}\ (\frac{ab}p)=1,
\\p^{(p-3)/4}&\text{if}\ (\frac{ab}p)=-1,\end{cases}$$
where $(\f{\cdot}p)$ denotes the Legendre symbol.
We also pose some conjectures for further research.
\end{abstract}
\maketitle

\section{Introduction}
\setcounter{lemma}{0}
\setcounter{theorem}{0}
\setcounter{corollary}{0}
\setcounter{remark}{0}
\setcounter{equation}{0}

Let $p$ be an odd prime, and let $(\f{\cdot}p)$ be the Legendre symbol.
Let $d$ be any integer. Sun \cite{S19} introduced the determinants
$$S(d,p)=\det\l[\l(\f{j^2+dk^2}p\r)\r]_{1\ls j,k\ls(p-1)/2}$$
and
$$T(d,p)=\det\l[\l(\f{j^2+dk^2}p\r)\r]_{0\ls j,k\ls(p-1)/2},
$$
and determined the Legendre symbols
$$\l(\f{S(d,p)}p\r)\ \t{and}\ \ \l(\f{T(d,p)}p\r).$$
Namely, the author \cite[Theorem 1.2]{S19} showed that
\begin{equation*}\l(\f{S(d,p)}p\r)=\begin{cases}(\f{-1}p)&\t{if}\ (\f dp)=1,
\\0&\t{if}\ (\f dp)=-1,\end{cases}
\end{equation*}
and
\begin{equation*}\l(\f{T(d,p)}p\r)=\begin{cases}(\f{2}p)&\t{if}\ (\f dp)=1,
\\1&\t{if}\ (\f dp)=-1.\end{cases}
\end{equation*}
D. Grinberg, the author and L. Zhao \cite{GSZ} proved that
if $p>3$ then
$$\det\l[(j^2+dk^2)\l(\f{j^2+dk^2}p\r)\r]_{0\ls j,k\ls(p-1)/2}\eq0\pmod p.$$
For any positive integer $n$ with $(p-1)/2\ls n\ls p-1$, we introduce the determinants
\begin{equation}S_n(d,p)=\det\l[(j^2+dk^2)^n\r]_{1\ls j,k\ls(p-1)/2}
\end{equation}
and
\begin{equation}T_n(d,p)=\det\l[(j^2+dk^2)^n\r]_{0\ls j,k\ls(p-1)/2}.
\end{equation}
Note that
$$S_{(p-1)/2}(d,p)\eq S(d,p)\pmod p,\ \ T_{(p-1)/2}(d,p)\eq T(d,p)\pmod p,$$
and
$$T_{(p+1)/2}(d,p)\eq\det\l[(j^2+dk^2)\l(\f{j^2+dk^2}p\r)\r]_{0\ls j,k\ls(p-1)/2}\pmod p.$$
When $p>3$ and $p\nmid d$, the author \cite[Conjecture 4.5(iii)]{S19} conjectured that
$$\l(\f{S_{(p+1)/2}(d,p)}p\r)=\begin{cases}(\f dp)^{(p-1)/4}&\t{if}\ p\eq1\pmod 4,
\\(\f dp)^{(p+1)/4}(-1)^{(h(-p)-1)/2}&\t{if}\ p\eq3\pmod4,
\end{cases}$$
where $h(-p)$ denotes the class number of the imaginary quadratic field $\Q(\sqrt{-p})$;
this was confirmed by H.-L. Wu, Y.-F. She and L.-Y. Wang \cite{WSW} in 2022.

\begin{theorem}\label{Th1.1} Let $p>3$ be a prime, and let $d\in\Z$.

{\rm (i)} Let $\bar S(d,p)$ be the determiant
obtained from $\det[(\f{j^2+dk^2}p)]_{1\ls j,k\ls(p-1)/2}$ by replacing all the entries in the first row by $1$. If $(\f dp)=1$, then
$$\bar S(d,p)=-S(d,p).$$ When $(\f dp)=-1$, we have
\begin{equation}\label{d/p}\bar S(d,p)=\f2{p-1}T(d,p)=\f{p-1}2\det\l[\l(\f{j^2+dk^2}p\r)\r]_{2\ls j,k\ls(p-1)/2}.
\end{equation}

{\rm (ii)}  For any integer $n$ with $(p-1)/2<n<p-1$, we have
\begin{equation} T_n(d,p)\eq0\pmod p.
\end{equation}
\end{theorem}
\begin{remark}
Part (ii) of Theorem \ref{Th1.1} extends \cite[Theorem 1.1]{GSZ}.
\end{remark}

For any prime $p\eq3\pmod4$, Sun \cite{S19} proved that
$$S_{p-2}(1,p)\eq\det\l[\f1{j^2+k^2}\r]_{1\ls j,k\ls(p-1)/2}\eq\l(\f 2p\r)\pmod p.$$
In contrast with this, we get the following result.

\begin{theorem} \label{Th-ij} Let $p$ be an odd prime, and let $d\in\Z$ with $(\f{-d}p)=-1$.

{\rm (i)} We have
\begin{equation}\label{p-2}\l(\f{S_{p-2}(d,p)}p\r)=\l(\f 2p\r).
\end{equation}
Moreover,
\begin{equation}\label{dk^2}\det\l[\f1{j^2+dk^2}\r]_{1\ls j,k\ls(p-1)/2}
\eq\begin{cases}d^{(p-1)/4}\pmod{p}&\t{if}\ p\eq1\pmod4,\\(-1)^{(p+1)/4}\pmod p&\t{if}\ p\eq3\pmod4.
\end{cases}
\end{equation}

{\rm (ii)} We have
\begin{equation}\label{Sp-3=1}\l(\f{S_{p-3}(d,p)}p\r)=\f{1-(\f{-1}p)}2.\end{equation}
Moreover, when $p\eq3\pmod4$ we have
\begin{equation}\label{i2+j2}\det\l[\f1{(j^2+dk^2)^2}\r]_{1\ls j,k\ls (p-1)/2}\eq\f14\prod_{r=1}^{(p-3)/4}\l(r+\f14\r)^2\pmod p.
\end{equation}
\end{theorem}

Let $p$ be an odd prime, and let $a,b\in\Z$ with $p\nmid ab$. The author \cite{RJ} introduced
\begin{equation}T^{(0)}_p(a,b,x)=\det\l[x+\tan\pi\f{aj^2+bk^2}p\r]_{0\ls j,k\ls (p-1)/2}
\end{equation}
and
\begin{equation}T^{(1)}_p(a,b,x)=\det\l[x+\tan\pi\f{aj^2+bk^2}p\r]_{1\ls j,k\ls (p-1)/2},
\end{equation}
and simply denote $T^{(0)}_p(a,b,0)$ and $T^{(1)}_p(a,b,0)$ by $T^{(0)}_p(a,b)$ and $T^{(1)}_p(a,b)$, respectively.
When $p>3$ and $p\eq3\pmod4$, the author \cite[Theorem 1.1(ii)]{RJ} proved that
\begin{equation}\label{Tp03}
T_p^{(0)}(a,b,x)=\begin{cases}2^{(p-1)/2}p^{(p+1)/4}&\t{if}\ (\f{ab}p)=1,
\\p^{(p+1)/4}&\t{if}\ (\f{ab}p)=-1.\end{cases}
\end{equation}
When $p\eq1\pmod4$, by \cite[Theorem 1.1(i)]{RJ} we have
 \begin{equation}\label{Tp+-}\begin{aligned}T_p^{(1)}(a,b,x)&=T_p^{(1)}(a,b)
 \\&=\begin{cases}\l(\f{2c}p\r)p^{(p-3)/4}\ve_p^{(\f ap)(2-(\f2p))h(p)}&\t{if}\ p\mid b-ac^2\ \t{with}\ c\in\Z,
 \\\pm2^{(p-1)/2}p^{(p-3)/4}&\t{if}\ (\f{ab}p)=-1,\end{cases}
 \end{aligned}
 \end{equation}
 where $\ve_p$ and $h(p)$ are the fundamental unit and the class number of the real quadratic field
 $\Q(\sqrt p)$, respectively.
As a supplement to \cite[Theorem 1.1]{RJ}, we obtain the following result.

\begin{theorem}\label{Th1.4} Let $p>3$ be a prime, and let $a,b\in\Z$ with $p\nmid ab$.

{\rm (i)} Assume that $p\eq1\pmod4$.
If $(\f{ab}p)=1$ and $ac^2\eq b$ with $c\in\Z$, then
\begin{equation}T_p^{(0)}(a,b,x)=\l(\f{2c}p\r)p^{(p+1)/4}\ve_p^{(\f ap)((\f 2p)-2)h(p)}x.
\end{equation}
 If $(\f{ab}p)=-1$, then
 \begin{equation}\label{T-da}T_p^{(1)}(a,b)=-\da(ab,p)2^{(p-1)/2}p^{(p-3)/4}
 \end{equation}
 and
 \begin{equation}T_p^{(0)}(a,b,x)=pT_p^{(1)}(a,b)x=-\da(ab,p)2^{(p-1)/2}p^{(p+1)/4}x,
 \end{equation}
 where
 \begin{equation}\label{delta}\da(c,p)=\begin{cases}1&\t{if}\ c^{(p-1)/4}\eq\f{p-1}2!\pmod p,\\-1&\t{otherwise}.
 \end{cases}\end{equation}

 {\rm (ii)} Suppose that $p\eq3\pmod4$.
 Then
 \begin{equation}\label{ConjT}T_p^{(1)}(a,b,x)
=\begin{cases}-2^{(p-1)/2}p^{(p-3)/4}x&\t{if}\ (\f{ab}p)=1,
\\p^{(p-3)/4}x&\t{if}\ (\f{ab}p)=-1.\end{cases}
\end{equation}
\end{theorem}
\begin{remark} In light of Theorem \ref{Th1.4} and \cite[Theorem 1.1]{RJ}, for any prime $p>3$ and
$a,b\in\Z$ with $p\nmid ab$, we have completed determined the exact values of $T_p^{(0)}(a,b,x)$
and $T_p^{(1)}(a,b,x)$.
\end{remark}

Let $p>3$ be a prime, and let $a,b\in\Z$ with $(\f{-ab}p)=-1$.
Define
\begin{equation}C_p(a,b,x)=\det\l[x+\cot\pi\f{aj^2+bk^2}p\r]_{1\ls j,k\ls(p-1)/2}.
\end{equation}
By \cite[Theorem 1.3]{RJ},
\begin{equation}\label{C}C_p(a,b,x)=\begin{cases} T_p^{(1)}(a,b)/(-p)^{(p-1)/4}=\pm 2^{(p-1)/2}/\sqrt p&\t{if}\ p\eq1\pmod4,
\\(-1)^{(h(-p)+1)/2}(\f ap)2^{(p-1)/2}/\sqrt p&\t{if}\ p\eq3\pmod4,
\end{cases}\end{equation}
In the case $p\eq1\pmod 4$, with the aid of \eqref{T-da} we have
\begin{equation}\label{cot}C_p(a,b,x)=(-1)^{(p+3)/4}\da(ab,p)\f{2^{(p-1)/2}}{\sqrt p}.
\end{equation}

Now we state our last two theorems.

\begin{theorem} \label{Th-tan} Let $p>3$ be a prime, and let $a,b\in\Z$ with $p\nmid ab$.
Let $\bar T_p(a,b,x)$ denote the determinant obtained from
$$T_p^{(0)}(a,b,x)=\det\l[x+\tan\pi\f{aj^2+bk^2}p\r]_{0\ls j,k\ls(p-1)/2}$$
 via replacing all the entries in the first row by $1$.

{\rm (i)} Suppose that $p\eq1\pmod 4$. If $(\f{ab}p)=1$ and $ac^2\eq b\pmod p$ with $c\in\Z$, then
\begin{equation}\bar T_p(a,b,x)=\l(\f{2c}p\r)p^{(p-1)/4}.
\end{equation}
If $(\f{ab}p)=-1$, then
\begin{equation}\bar T_p(a,b,x)=-\da(ab,p)2^{(p-1)/2}p^{(p-1)/4}\ve_p^{(\f{2a}p)h(p)}.
\end{equation}

{\rm (ii)} When $p\eq3\pmod4$, we have
\begin{equation}\label{Tp3}\bar T_p(a,b,x)=(-1)^{\f{p+1}4+\f{h(-p)+1}2}\l(\f ap\r)2^{(1+(\f{ab}p))\f{p-1}4}p^{(p-1)/4}.
\end{equation}
\end{theorem}

\begin{theorem} \label{Th-cot} Let $p>3$ be a prime, and let $a,b\in\Z$ with $(\f{-ab}p)=-1$.
Let $\bar C_p(a,b,x)$ denote the determinant of the matrix $[c_{jk}]_{0\ls j,k\ls(p-1)/2}$,
where
$$ c_{jk}=\begin{cases}1&\t{if}\ j=0,\\ x+\cot\pi(aj^2+bk^2)/p&\t{if}\ j>0.\end{cases}$$
Then
\begin{equation}\bar C_p(a,b,x)=\f{2^{(p-1)/2}}{\sqrt p}\times
\begin{cases}(-1)^{(p+3)/4}\da(ab,p)\ve_p^{(\f ap)2h(p)}
&\t{if}\ p\eq1\pmod4,\\(-1)^{(h(-p)-1)/2}(\f ap)&\t{if}\ p\eq3\pmod4.\end{cases}
\end{equation}
\end{theorem}

We are going to prove Theorems 1.1-1.2 in the next section.
Based on two auxiliary theorems in Section 3, we will prove Theorem 1.3 in Section 4.
Our proofs of Theorems 1.4-1.5 will be given in Section 5.
In Section 6 we pose several conjectures on determinants for further research.

\section{Proofs of Theorems 1.1-1.2}
\setcounter{lemma}{0}
\setcounter{theorem}{0}
\setcounter{corollary}{0}
\setcounter{remark}{0}
\setcounter{equation}{0}

We need the following known lemma (cf. \cite[p.\,58]{BEW}).

\begin{lemma}\label{Lem-abc} Let $p$ be an odd prime, and let $a,b,c\in\Z$ with $p\nmid a$. Then
\begin{equation}\sum_{x=0}^{p-1}\l(\f{ax^2+bx+c}p\r)=\begin{cases}(p-1)(\f ap)&\t{if}\ p\mid b^2-4ac,
\\-(\f ap)&\t{if}\ p\nmid b^2-4ac.\end{cases}\end{equation}
\end{lemma}

\medskip
\noindent{\tt Proof of Theorem 1.1(i)}. By Lemma \ref{Lem-abc}, for each $k=1,\ldots,(p-1)/2$ we have
\begin{equation}\label{jk^2}\sum_{j=1}^{(p-1)/2}\l(\f{j^2+dk^2}p\r)=\f12\(\sum_{j=0}^{p-1}\l(\f{j^2+dk^2}p\r)-\l(\f{dk^2}p\r)\)=-\f{1+(\f dp)}2.
\end{equation}
Thus, for the determinant $T(d,p)=|(\f{j^2+dk^2}p)|_{0\ls j,k\ls(p-1)/2}$, if we add all the rows below the second row to the second row, then the second row becomes
$$\l(\f{p-1}2,-\f{1+(\f dp)}2,\ldots,-\f{1+(\f dp)}2\r)$$
while the first row is
$$\l(0,\l(\f dp\r),\ldots,\l(\f dp\r)\r).$$
Therefore, in the case $(\f dp)=-1$, we have
$$T(d,p)=\f{p-1}2\bar S(d,p).$$

Now we consider the case $(\f dp)=1$. If we add to the second row of $T(d,p)$ all the other rows,
then the second row becomes $(\f{p-1}2,0,\ldots,0)$ by \eqref{jk^2} while the first row is
$(0,1,\ldots,1)$. It follows that
$$T(d,p)=-\f{p-1}2\bar S(d,p).$$
By \cite[(1.20)]{S19},
$$T(d,p)=\f{p-1}2S(d,p).$$
Combining the last two equalities, we get $\bar S(d,p)=-S(d,p)$.

By Lemma \ref{Lem-abc}, for any $j=1,\ldots,(p-1)/2$ we have
\begin{equation}
\label{sum-k} \sum_{k=1}^{(p-1)/2}\l(\f{j^2+dk^2}p\r)=\f12\(\sum_{k=0}^{p-1}\l(\f{j^2+dk^2}p\r)-1\)=-\f{(\f dp)+1}2.
\end{equation}
Suppose that $(\f dp)=-1$. If we add to the first column
of $\bar S(s,p)$ all the other columns, then the first column turns out to be $(\f{p-1}2,0,\ldots,0)^T$
by \eqref{sum-k}.
Therefore,
$$\bar S(d,p)=\f{p-1}2\det\l[\l(\f{j^2+dk^2}p\r)\r]_{2\ls j,k\ls(p-1)/2}.$$

Combining the above, we have completed our proof of Theorem 1.1(i). \qed

\medskip
\noindent{\tt Proof of Theorem 1.1(ii)}. Let $k\in\{1,\ldots,(p-1)/2\}$.
In view of the binomial theorem, we have
\begin{align*}\sum_{j=1}^{(p-1)/2}(j^2+dk^2)^n&=\sum_{j=1}^{(p-1)/2}\sum_{r=0}^n\bi nrj^{2r}(dk^2)^{n-r}
\\&\eq\sum_{r=0}^n\bi nr(dk^2)^{n-r}\f12\sum_{j=1}^{(p-1)/2}\l(j^{2r}+(p-j)^{2r}\r)
\\&=\sum_{r=0}^n\bi nr(dk^2)^{n-r}\sum_{j=1}^{p-1}j^{2r}\pmod p.
\end{align*}
By a well known result (cf. \cite[Section 15.2, Lemma 2]{IR}),
$$\sum_{j=1}^{p-1}j^{2r}\eq\begin{cases}-1\pmod p&\t{if}\ p-1\mid 2r,
\\0\pmod p&\t{otherwise}.\end{cases}$$
As $(p-1)/2<n<p-1$,  for $r\in\{0,\ldots,n\}$ we have
$$p-1\mid 2r\iff \f{p-1}2\ \big|\ r\iff r=0\ \t{or}\ r=\f{p-1}2.$$
Thus
\begin{align*}\sum_{j=1}^{(p-1)/2}(j^2+dk^2)^n&\eq-\sum_{r\in\{0,(p-1)/2\}}\bi nr(dk^2)^{n-r}
\\=&-(dk^2)^n-\bi n{(p-1)/2}(dk^2)^{n-(p-1)/2}\pmod p
\end{align*}
and hence
$$\l(1+\l(\f dp\r)\bi n{(p-1)/2}\r)(dk^2)^n+\sum_{j=1}^{(p-1)/2}(j^2+dk^2)^n\eq0\pmod p.
$$

As $p-1\nmid 2n$, we also have
$$\l(1+\l(\f dp\r)\bi n{(p-1)/2}\r)(d0^2)^n+\sum_{j=1}^{(p-1)/2}(j^2+d0^2)^n\eq\f12\sum_{j=1}^{p-1}j^{2n}\eq0\pmod p.
$$
Combining this with the last paragraph, we see that
$$\l(1+\l(\f dp\r)\bi n{(p-1)/2}\r)t_{0k}+\sum_{j=1}^{(p-1)/2}t_{jk}\eq0\pmod p
$$
for all $k=0,\ldots,(p-1)/2$, where $t_{jk}=(j^2+dk^2)^n$. Therefore
$$T_n(d,p)=\det[t_{jk}]_{0\ls j,k\ls(p-1)/2}\eq0\pmod p$$
as desired. \qed

The following well known result can be found in the survey \cite[(5.5)]{K}.

\begin{lemma}[Cauchy] \label{Lem-C} We have
\begin{equation} \det\l[\f1{x_j+y_k}\r]_{1\ls j,k\ls n}=\f{\prod_{1\ls j<k\ls n}(x_j-x_k)(y_j-y_k)}{\prod_{j=1}^n\prod_{k=1}^n(x_j+y_k)}.
\end{equation}
\end{lemma}

Let $p$ be an odd prime. In view of Wilson's theorem,
$$\prod_{k=1}^{(p-1)/2}k(p-k)=(p-1)!\eq-1\pmod p$$
and hence
\begin{equation}\label{p/2}\l(\f{p-1}2!\r)^2\eq(-1)^{(p+1)/2}\pmod p.
\end{equation}
By \cite[(1.5)]{S19}, we have
\begin{equation}\label{k^2-j^2}\prod_{1\ls j<k\ls (p-1)/2}(k^2-j^2)\eq\begin{cases}-\f{p-1}2!\pmod p&\t{if}\ p\eq1\pmod4,
\\1\pmod p&\t{if}\ p\eq3\pmod4.\end{cases}
\end{equation}
Therefore
\begin{equation}\label{kj^2}\prod_{1\ls j<k\ls (p-1)/2}(k^2-j^2)^2\eq(-1)^{(p+1)/2}\pmod p.
\end{equation}

\medskip
\noindent {\tt Proof of Theorem \ref{Th-ij}(i)}.
Let $n=(p-1)/2$. By Lemma \ref{Lem-C} and \eqref{kj^2}, we have
\begin{align*}\det\l[\f1{j^2+dk^2}\r]_{1\ls j,k\ls n}&=\f{\prod_{1\ls j<k\ls n}(k^2-j^2)(dk^2-dj^2)}{\prod_{j=1}^n\prod_{k=1}^n(j^2+dk^2)}
\\&=\f{d^{n(n-1)/2}}{\Pi}\prod_{1\ls j<k\ls n}(k^2-j^2)^2
\\&\eq(-1)^{n+1}\f{d^{n(n-1)/2}}{\Pi}\pmod p,
\end{align*}
where
$$\Pi:=\prod_{k=1}^n\(k^{2n}\prod_{j=1}^n\l(\f {j^2}{k^2}+d\r)\)\eq\prod_{k=1}^n\prod_{x=1}^n(x^2+d)\pmod p.$$

Note that
$$\prod_{x=1}^n(x^2+d)\eq(-1)^{n+1}2\pmod p$$
by \cite[Lemma 3.1]{S19}. Thus
$$\Pi\eq((-1)^{n+1}2)^n\eq 2^n\eq\l(\f 2p\r)=(-1)^{(p^2-1)/8}\pmod p.$$
If $p\eq1\pmod 4$, then $2\mid n$ and hence
$$d^{n(n-1)/2}=(d^n)^{n/2-1}d^{n/2}\eq\l(\f dn\r)^{n/2-1}d^{n/2}=(-1)^{n/2-1}d^{(p-1)/4}\pmod p.$$
If $p\eq3\pmod 4$, then $2\nmid n$ and hence
$$d^{n(n-1)/2}=(d^n)^{(n-1)/2}\eq\l(\f dp\r)^{(n-1)/2}=1\pmod p.$$
Therefore
$$\f{d^{n(n-1)/2}}{\Pi}\eq\begin{cases}-d^{(p-1)/4}\pmod p&\t{if}\ p\eq1\pmod4,
\\(-1)^{(p+1)/4}\pmod p&\t{if}\ p\eq3\pmod4.\end{cases}$$
Combining this with the first paragraph in the proof, we immediately obtain
the congruence \eqref{dk^2}, which clearly implies \eqref{p-2}.
This concludes the proof. \qed

Recall that the permanent of an $n\times n$ matrix $A=(a_{j,k})_{1\ls j,k\ls n}$ over a field is given by
$$\per(A)=\per[a_{j,k}]_{1\ls j,k\ls n}=\sum_{\sigma\in S_n}\prod_{j=1}^na_{j,\sigma(j)}.$$

\begin{lemma}\label{per-d} Let $p$ be an odd prime, and let $d\in\Z$ with $(\f {-d}p)=-1$. Then
\begin{equation}\label{per}\begin{aligned}&\ \per\l[\f1{j^2+dk^2}\r]_{1\ls j,k\ls(p-1)/2}
\\\eq&\ \begin{cases}0\pmod p&\t{if}\ p\eq1\pmod 4,
\\(-1)^{(p-3)/4}\prod_{r=1}^{(p-3)/4}(r+\f14)^2\pmod p&\t{if}\ p\eq3\pmod4.
\end{cases}
\end{aligned}\end{equation}
\end{lemma}
\Proof. Let $g$ be a primitive root modulo $p$, and set $n=(p-1)/2$. Then
those $g^{2k}\ (k=1,\ldots,n)$ are incongruent quadratic residues modulo $p$.
Thus
\begin{align*} \per\l[\f1{j^2+dk^2}\r]_{1\ls j,k\ls n}
&=\f1{\prod_{k=1}^n k^2}\per\l[\f1{1+dk^2/j^2}\r]_{1\ls j,k\ls n}
\\&\eq\f1{(n!)^2}\per\l[\f1{1+dg^{2(j-k)}}\r]_{1\ls j,k\ls n}
\\&\eq(-1)^{n-1}\prod_{r=1}^n\l(\f{n(-d)^n}{1-(-d)^n}+r\r)\pmod p
\end{align*}
by \eqref{p/2} and \cite[Theorem 1.3(i)]{S-perm}.
As $(-d)^n\eq(\f{-d}p)=-1\pmod p$, from the above we get
\begin{equation}\label{1/4}\per\l[\f1{j^2+dk^2}\r]_{1\ls j,k\ls n}
\eq(-1)^{n-1}\prod_{r=1}^n\l(r+\f14\r)\pmod p.
\end{equation}
If $p\eq1\pmod 4$, then $r+1/4\eq0\pmod p$ for $r=(p-1)/4$.
When $p\eq3\pmod4$, we have
\begin{align*}\prod_{r=1}^n\l(r+\f14\r)&=\l(\f{p-1}2+\f14\r)\prod_{r=1}^{(p-3)/4}\l(r+\f14\r)\l(\f{p-1}2-r+\f14\r)
\\&\eq \f{(-1)^{(p+1)/4}}4\prod_{r=1}^{(p-3)/4}\l(r+\f14\r)^2\pmod p.
\end{align*}
Therefore \eqref{1/4} implies the desired congruence \eqref{per}. \qed

The following result due to  Borchardt can be found in \cite{Sin}.

\begin{lemma} \label{Lem-B} We have
\begin{equation}\label{Bor}
\det\l[\f1{(x_j+y_k)^2}\r]_{1\ls j,k\ls n}=\det\l[\f1{x_j+y_k}\r]_{1\ls j,k\ls n}\per\l[\f1{x_j+y_k}\r]_{1\ls j,k\ls n}.
\end{equation}
\end{lemma}

\medskip
\noindent{\tt Proof of Theorem} \ref{Th-ij}(ii). Combining \eqref{dk^2},
and Lemmas \ref{per-d} and \ref{Lem-B}, we immediately get the desired results. \qed

\section{Two Auxiliary Theorems}
\setcounter{lemma}{0}
\setcounter{theorem}{0}
\setcounter{corollary}{0}
\setcounter{remark}{0}
\setcounter{equation}{0}

Our first auxiliary theorem is as follows.

\begin{theorem}\label{Th-m} Let $p$ be an odd prime, and let $k,m\in\Z^+=\{1,2,3,\ldots\}$ with $km=p-1$.
Let $G$  be the multiplicative group $\{r+p\Z:\ r=1,\ldots,p-1\}$
and let $H$ be its subgroup $\{x^m+p\Z:\ x=1,\ldots,p-1\}$ of order $k$. Suppose that all the $m$ distinct cosets of $H$ in $G$
are
$$\{a_{1j}+p\Z:\ j=1,\ldots,k\},\ \ldots,\ \{a_{mj}+p\Z:\ j=1,\ldots,k\}$$
with $1\ls a_{i1}<\ldots<a_{ik}\ls p-1$ for all $i=1,\ldots,m$. Then
\begin{equation}\label{m}\begin{aligned}&\prod_{i=1}^m\prod_{1\ls s<t\ls k}(a_{it}-a_{is})
\\\eq&\begin{cases}(-1)^{\f{p+1}2\cdot\f{p-1}{2m}+\lfloor\f{p-3}4\rfloor}\,\f{p-1}2!\pmod p&\t{if}\ p\eq1\pmod{2m},\\(-1)^{\f{p+1}2\cdot\f{p-1-m}{2m}}\pmod p&\t{if}\ p\eq 1+m\pmod{2m}.
\end{cases}
\end{aligned}
\end{equation}
\end{theorem}
\Proof. Set
$$R_m=\{1\ls r\ls p-1:\ x^m\eq r\pmod p\ \t{for some}\ x=1,\ldots,p-1\}.$$
Then $H=\{r+p\Z:\ r\in R_m\}$ and $|H|=|R_m|=(p-1)/m=k$.
Note that
$$\prod_{i=1}^m\prod_{1\ls s<t\ls k}(a_{it}-a_{is})=\prod_{d=1}^{p-1}d^{e_d},$$
where
\begin{align*}e_d:&= |\{1\ls x<p-d:\ \{x,x+d\}\se\{a_{i1},\ldots,a_{ik}\}\ \t{for some}\ i=1,\ldots,m\}|
\\&=\l|\l\{1\ls x<p-d:\ \f{x+d}x\eq r\pmod p\ \t{for some}\ r\in R_m\r\}\r|.
\end{align*}
Clearly,
$$\prod_{d=1}^{p-1}d^{e_d}=\prod_{d=1}^{(p-1)/2}d^{e_d}(p-d)^{e_{p-d}}\eq(-1)^{\sum_{d=1}^{(p-1)/2}e_{p-d}}
\prod_{d=1}^{(p-1)/2}d^{e_d+e_{p-d}}\pmod p.$$

For any $d\in\{1,\ldots,p-1\}$, obviously
\begin{align*}e_{p-d}&=\l|\l\{1\ls x<d:\ \f{x+p-d}x\eq r\pmod p\ \t{for some}\ r\in R_m\r\}\r|
\\&=\l|\l\{p-d<y<p:\ 1+\f{p-d}{p-y}\eq r\pmod p\ \t{for some}\ r\in R_m\r\}\r|
\\&=\l|\l\{p-d\ls y<p:\ \f{y+d}{y}\eq r\pmod p\ \t{for some}\ r\in R_m\r\}\r|
\end{align*}
and hence
\begin{align*}e_d+e_{p-d}&=\l|\l\{1\ls x<p:\ 1+\f dx\eq r\pmod p\ \t{for some}\ r\in R_m\r\}\r|
\\&=|\{1<y<p:\ y\eq r\pmod p\ \t{for some}\ r\in R_m\}|
\\&=|R_m|-1=k-1.
\end{align*}

Observe that
$$\sum_{d=1}^{(p-1)/2}e_{p-d}
=\sum_{d=1}^{(p-1)/2}\l|\l\{1\ls x<d:\ \f {d-x}x\eq-r\pmod p\ \t{for some}\ r\in R_m\r\}\r|$$
coincides with
$$\l|\l\{(x,y)\in(\Z^+)^2:\ x+y\ls\f{p-1}2\ \t{and}\ \f yx\eq-r\pmod p\ \ \t{for some}\ r\in R_m\r\}\r|.
$$
As $H$ is a multiplicative group, given $x,y\in\{1,\ldots,p-1\}$ we have
$$\f yx\eq-r\pmod p\ \ \t{for some}\ r\in R_m
\iff \f xy\eq-r\pmod p\ \ \t{for some}\ r\in R_m.$$
Therefore,
$\sum_{d=1}^{(p-1)/2}e_{p-d}$
has the same parity with
\begin{align*}&\ \l|\l\{x\in\Z^+:\ x+x\ls\f{p-1}2\ \t{and}\ \f xx\eq-r\pmod p\ \ \t{for some}\ r\in R_m\r\}\r|
\\=&\ \l|\l\{1\ls x<\f p4:\ p-1\in R_m\r\}\r|
=\l|\l\{1\ls x<\f p4:\ (-1)^{(p-1)/m}=1\r\}\r|
\\=&\begin{cases}\lfloor(p-1)/4\rfloor&\t{if}\ 2\mid k,\\0&\t{if}\ 2\nmid k,\end{cases}
\end{align*}
and hence
$$(-1)^{\sum_{d=1}^{(p-1)/2}e_{p-d}}=(-1)^{(k-1)\lfloor\f{p-1}4\rfloor}.$$

Combining the above, we see that
$$\prod_{i=1}^m\prod_{1\ls s<t\ls k}(a_{it}-a_{is})\eq(-1)^{(k-1)\lfloor\f{p-1}4\rfloor}\prod_{d=1}^{(p-1)/2}d^{k-1}\pmod p.$$
Recall that
$$\l(\f{p-1}2!\r)^2\eq(-1)^{(p+1)/2}\pmod p$$
by Wilson's theorem. So, by the last two congruences we immediately obtain the desired congruence \eqref{m}. \qed

Theorem \ref{Th-m} in the case $m=2$ yields the following result.

\begin{corollary}\label{Lem-RN} Let $p=2n+1$ be an odd prime, and write $$\{1,\ldots,p-1\}=\{a_1,\ldots,a_n\}\cup\{b_1,\ldots,b_n\}$$
with $a_1<\ldots<a_n$ and $b_1<\ldots<b_n$ such that $a_1,\ldots,a_n$
are quadratic residues modulo $p$, and $b_1,\ldots,b_n$ are quadratic nonresidues modulo $p$.
Then
\begin{equation}\label{ab-cong} \prod_{1\ls j<k\ls n}(a_k-a_j)(b_k-b_j)\eq \begin{cases}-n!\pmod p&\t{if}\ p\eq1\pmod4,
\\1\pmod p&\t{if}\ p\eq3\pmod4.\end{cases}
\end{equation}
\end{corollary}

For any odd prime $p$ and integer $a\not\eq0\pmod p$, we define
$$s_p(a)=(-1)^{|\{\{j,k\}:\ 1\ls j<k\ls (p-1)/2\ \t{and}\ \{aj^2\}_p>\{ak^2\}_p\}|},$$
where $\{m\}_p$ denotes the least nonnegative residue of an integer $m$ modulo $p$.
The author \cite[Theorem 1.4(i)]{S19p} determined $s_p(1)$ in the case $p\eq3\pmod4$.
When $p\eq1\pmod4$, H.-L. Wu \cite{Wu} deduced a complicated formula for $s_p(1)$ modulo $p$,
which involves the fundamental unit $\ve_p$ and the class numbers of the quadratic fields $\Q(\sqrt{\pm p})$.

Based on Corollary \ref{Lem-RN}, we get the following result.

\begin{lemma} Let $p$ be an odd prime, and let $a,b\in\Z$ with $(\f ap)=1$ and $(\f bp)=-1$.
Then
\begin{equation}s_p(a)s_p(b)=\begin{cases}(-1)^{(p+3)/4}\da(ab,p)&\t{if}\ p\eq1\pmod4,
\\(-1)^{(p-3)/4}&\t{if}\ p\eq3\pmod4.\end{cases}
\end{equation}
\end{lemma}
\Proof. Let $n=(p-1)/2$, and write $\{1,\ldots,p-1\}=\{a_1,\ldots,a_n\}\cup\{b_1,\ldots,b_n\}$ with $a_1<\ldots<a_n$ and $b_1<\ldots<b_n$ such that $a_1,\ldots,a_n$ are quadratic residues modulo $p$
and $b_1,\ldots,b_n$ are quadratic nonresidues modulo $p$. As
$$\{\{aj^2\}_p:\ j=1,\ldots,n\}=\{a_1,\ldots,a_n\}$$
and
$$\{\{bj^2\}_p:\ j=1,\ldots,n\}=\{b_1,\ldots,b_n\},$$
we have
$$s_p(a)s_p(b)=\prod_{1\ls j<k\ls n}\f{\{ak^2\}_p-\{aj^2\}_p}{a_k-a_j}\times\prod_{1\ls j<k\ls n}\f{\{bk^2\}_p-\{aj^2\}_p}{b_k-b_j}$$
and hence
\begin{align*}&\ s_p(a)s_p(b)\prod_{1\ls j<k\ls n}(a_k-a_j)(b_k-b_j)
\\\eq&\ \prod_{1\ls j<k\ls n}(ak^2-aj^2)(bk^2-bj^2)=(ab)^{n(n-1)/2}\prod_{1\ls j<k\ls n}(k^2-j^2)^2
\pmod p.
\end{align*}
Note that
$$(ab)^n\eq\l(\f{ab}p\r)=-1\pmod p.$$
By \cite[(1.5)]{S19} we have
\begin{equation}\label{k^2-j^2}\prod_{1\ls j<k\ls n}(k^2-j^2)\eq\begin{cases}-n!\pmod p&\t{if}\ p\eq1\pmod4,
\\1\pmod p&\t{if}\ p\eq3\pmod4,\end{cases}
\end{equation}
and hence
$$\prod_{1\ls j<k\ls n}(k^2-j^2)^2\eq(-1)^{(p+1)/2}\pmod p.$$
Therefore
\begin{align*}&\ s_p(a)s_p(b)\prod_{1\ls j<k\ls n}(a_k-a_j)(b_k-b_j)
\\\eq&\begin{cases}(-1)^{n/2-1}(ab)^{n/2}\times(-1)\pmod p&\t{if}\ p\eq1\pmod4,
\\(-1)^{(n-1)/2}\times1\pmod p&\t{if}\ p\eq3\pmod4.
\end{cases}
\end{align*}
Combining this with \eqref{ab-cong}, we obtain that
\begin{align*} s_p(a)s_p(b)
\eq\begin{cases} (-1)^{n/2}(ab)^{n/2}/(-n!)\pmod p&\t{if}\ p\eq1\pmod4,\\ (-1)^{(n-1)/2}\pmod p&\t{if}\ p\eq3\pmod4.
\end{cases}\end{align*}
In the case $p\eq1\pmod 4$, we have
$$(ab)^n\eq-1\eq(n!)^2\pmod p$$
and hence
$$(ab)^{n/2}\eq \pm n!\pmod p,$$
therefore
$$s_p(a)s_p(b)=(-1)^{n/2+1}\da(ab,p)=(-1)^{(p+3)/4}\da(ab,p).$$
This concludes our proof. \qed

Now we are ready to present another auxiliary theorem.

\begin{theorem}\label{Au-Th} Let $p$ be a prime with $p\eq1\pmod4$, and let $\zeta=e^{2\pi i/p}$.
Let $a,b\in\Z$ with $(\f {ab}p)=-1$. Then, we have
\begin{equation}\label{zeta-ab}\prod_{1\ls j<k\ls {(p-1)/2}}(\zeta^{aj^2}-\zeta^{ak^2})(\zeta^{bj^2}-\zeta^{bk^2})
=-\da(ab,p)p^{(p-3)/4}
\end{equation}
and
\begin{equation}\label{cot-ab}\begin{aligned}&\prod_{1\ls j<k\ls (p-1)/2}\l(\cot\pi\f{aj^2}p-\cot\pi\f{ak^2}p\r)\l(\cot\pi\f{bj^2}p-\cot\pi\f{bk^2}p\r)
\\&\qquad =\da(ab,p)(-1)^{(p+3)/4}\l(\f{2^{p-1}}p\r)^{(p-3)/4}.
\end{aligned}
\end{equation}
\end{theorem}
\begin{remark} For any prime $p>3$ with $p\eq3\pmod4$ and integer $a\not\eq0\pmod p$, the author \cite[part (ii) of Theorems 1.3-1.4]{S19} obtained closed forms for the products
$$\prod_{1\ls j<k\ls (p-1)/2}\l(e^{2\pi iaj^2/p}-e^{2\pi iak^2/p}\r)
\ \t{and}\ \prod_{1\ls j<k\ls(p-1)/2}\l(\cot\pi\f{aj^2}p-\cot\pi\f{ak^2}p\r).$$
\end{remark}

\medskip
\noindent{\tt Proof of Theorem} \ref{Au-Th}. Set $n=(p-1)/2$ and $\zeta=e^{2\pi i/p}$. By \cite[(4.2) and (4.3)]{S19p}, we have
$$\prod_{1\ls j<k\ls n}\sin\pi\f{a(k^2-j^2)}p=(-1)^{a(n+1)n/2}\l(\f i2\r)^{n(n-1)/2}\prod_{1\ls j<k\ls n}
(\zeta^{aj^2}-\zeta^{ak^2})$$
and
$$\f{\prod_{1\ls j<k\ls n}\sin\pi\f{a(k^2-j^2)}p}{\prod_{1\ls j<k\ls n}(\cot\pi\f{aj^2}p-\cot\pi\f{ak^2}p)}
=\l(\f p{2^{p-1}}\r)^{(n-1)/2}(-1)^{(a-1)n/2}\ve_p^{(\f ap)(1-n)h(p)}.$$
Therefore
\begin{equation}\label{cot-a}\prod_{1\ls j<k\ls n}\f{\cot\pi\f{aj^2}p-\cot\pi\f{ak^2}p}{\zeta^{aj^2}-\zeta^{ak^2}}
=\l(\f{2^n}p\r)^{(n-1)/2}i^{n(n+1)/2}\ve_p^{(\f ap)(1-n)h(p)}.
\end{equation}
Similarly,
\begin{equation}\label{cot-b}\prod_{1\ls j<k\ls n}\f{\cot\pi\f{bj^2}p-\cot\pi\f{bk^2}p}{\zeta^{bj^2}-\zeta^{bk^2}}
=\l(\f{2^n}p\r)^{(n-1)/2}i^{n(n+1)/2}\ve_p^{(\f bp)(1-n)h(p)}.
\end{equation}
Combining \eqref{cot-a} with \eqref{cot-b}, and noting $(\f ap)+(\f bp)=0$, we deduce that
\begin{equation}\prod_{1\ls j<k\ls n}\f{(\cot\pi \f{aj^2}p-\cot\pi\f{ak^2}p)(\cot\pi \f{bj^2}p-\cot\pi\f{bk^2}p)}{(\zeta^{aj^2}-\zeta^{ak^2})(\zeta^{bj^2}-\zeta^{bk^2})}=(-1)^{n/2}\l(\f{2^n}p\r)^{n-1}.
\end{equation}
So \eqref{zeta-ab} and \eqref{cot-ab} are equivalent.

By \cite[Theorem 1.3(i)]{S19p},
$$\prod_{1\ls j<k\ls n}(\zeta^{aj^2}-\zeta^{ak^2})=t_p(a)i^{n/2}p^{(n-1)/4}\ve_p^{(\f ap)\f{h(p)}2}$$
for some $t_p(a)\in\{\pm1\}$.
Combining this with \eqref{cot-a} we see that $t_p(a)$ coincides with the sign of the product
$$\prod_{1\ls j<k\ls n}\l(\cot\pi\f{aj^2}p-\cot\pi\f{ak^2}p\r)$$
which should be
$$(-1)^{|\{1\ls j<k\ls n:\ \{aj^2\}_p>\{ak^2\}_p\}|}=s_p(a).$$
Thus
$$\prod_{1\ls j<k\ls n}(\zeta^{aj^2}-\zeta^{ak^2})=s_p(a)i^{n/2}p^{(n-1)/4}\ve_p^{(\f ap)\f{h(p)}2}.$$
Similarly,
$$\prod_{1\ls j<k\ls n}(\zeta^{bj^2}-\zeta^{bk^2})=s_p(b)i^{n/2}p^{(n-1)/4}\ve_p^{(\f bp)\f{h(p)}2}.$$
Therefore
\begin{align*}&\ \prod_{1\ls j<k\ls n}(\zeta^{aj^2}-\zeta^{ak^2})(\zeta^{bj^2}-\zeta^{bk^2})
\\=&\ s_p(a)s_p(b)(-1)^{n/2}p^{(n-1)/2}=-\da(ab,p)p^{(p-3)/4}.
\end{align*}
This proves \eqref{zeta-ab}.

In view of the above, we have completed our proof of Theorem \ref{Au-Th}. \qed

\section{Proof of Theorem \ref{Th1.4}}
\setcounter{lemma}{0}
\setcounter{theorem}{0}
\setcounter{corollary}{0}
\setcounter{remark}{0}
\setcounter{equation}{0}

The following lemma is a known result (see, e.g., \cite[(1.12)]{S19p}).

\begin{lemma}\label{Lem-h} For any prime $p\eq1\pmod4$ and integer $a\not\eq0\pmod p$, we have
\begin{equation}\prod_{k=1}^{(p-1)/2}\l(1-e^{2\pi iak^2/p}\r)=\sqrt p\,\ve_p^{-(\f ap)h(p)}.
\end{equation}
\end{lemma}

\begin{lemma} \label{Lem3.1} Let $m,n\in\Z^+=\{1,2,3,\ldots\}$ with $2\nmid n$. Let $a_k,b_k\in\Z$ for $k=0,1,\ldots,m$. with $a_0+b_0=0$. Then
\begin{equation}\label{aj+bk}\begin{aligned}&\ \det\l[x+\tan\pi\f{a_j+b_k}n\r]_{0\ls j,k\ls m}-\det\l[\tan\pi\f{a_j+b_k}n\r]_{0\ls j,k\ls m}
\\=&\ x\det\l[\tan\pi\f{a_j+b_k}n\r]_{1\ls j,k\ls m}\times\prod_{k=1}^m\l(\tan\pi\f{a_k+b_0}n\times\tan\pi\f{a_0+b_k}n\r).
\end{aligned}
\end{equation}
\end{lemma}
\Proof. Let $a_{jk}=\tan\pi(a_j+b_k)/n$ for $j,k=0,\ldots,m$. By \cite[Lemma 2.1]{RJ}, we have
\begin{equation}\label{ab}\det[x+a_{jk}]_{0\ls j,k\ls m}-\det[a_{jk}]_{0\ls j,k\ls m}=x\det[b_{jk}]_{1\ls j,k\ls m},
\end{equation}
where $b_{jk}=a_{jk}-a_{j0}-a_{0k}+a_{00}$. Note that $a_{00}=\tan 0=0$ and recall the known identity
$$(1-\tan x_1\times\tan x_2)\tan(x_1+x_2)=\tan x_1+\tan x_2.$$
Then we have
\begin{align*}b_{jk}&=\tan \pi\f{a_j+b_k}n-\tan\pi\f{a_j+b_0}n-\tan\pi\f{a_0+b_k}n
\\&=\tan\pi\f{a_j+b_0}n\times\tan\pi\f{a_0+b_k}n\times\tan\pi\f{a_j+b_k}n.
\end{align*}
Thus
$$\det[b_{jk}]_{1\ls j,k\ls m}=\det\l[\tan\pi\f{a_j+b_k}n\r]_{1\ls j,k\ls m}\prod_{k=1}^m\l(\tan\pi\f{a_j+b_0}n\times\tan\pi\f{a_0+b_k}n\r).$$
Combining this with \eqref{ab}, we immediately obtain the desired identity \eqref{aj+bk}. \qed

 \medskip
 \noindent{\tt Proof of Theorem} 1.3(i). Let $n=(p-1)/2$, and let $a_{jk}=\tan \pi(aj^2+bk^2)/p$ for $j,k=0,n$.
 Set $q=n!$. By Wilson's theorem, we have $q^2\eq-1\pmod p$. Thus
 \begin{align*}T_p^{(0)}(a,b)&=\det\l[\tan\pi\f{a(qj)^2+b(qk)^2}p\r]_{0\ls j,k\ls n}
 \\&=\det\l[-\tan\pi\f{aj^2+bk^2}p\r]_{0\ls j,k\ls n}=-T_p^{(0)}(a,b)
 \end{align*}
 and hence $T_p^{(0)}(a,b)=0$ (which also follows from \cite[(1.3)]{RJ}).

 In view of the above and Lemma \ref{Lem3.1}, we have
 $$T_p^{(0)}(a,b,x)=xT_p^{(1)}(a,b)\prod_{k=1}^n\l(\tan\pi\f{ak^2}p\times\tan\pi\f{bk^2}p\r).$$
For any $x\in\Q$ with odd denominator, clearly
$$\tan\pi x=\f{2\sin \pi x}{2\cos\pi x}=\f{(e^{i\pi x}-e^{-i\pi x})/i}{e^{i\pi x}+e^{-i\pi x}}=i\f{1-e^{2\pi ix}}{1+e^{2\pi ix}}=i\f{(1-e^{2\pi ix})^2}{1-e^{2\pi i(2x)}}.$$
In view of this and
 Lemma \ref{Lem-h}, we deduce that
 \begin{align*}\prod_{k=1}^{n}\tan\pi\f{ak^2}p
 &=i^{n}\f{\prod_{k=1}^n(1-e^{2\pi i ak^2/p})^2}{\prod_{k=1}^n(1-e^{2\pi i(2a)k^2/p})}
 \\&==(i^2)^{n/2}\f{(\sqrt p\,\ve_p^{-(\f ap)h(p)})^2}{\sqrt p\,\ve_p^{-(\f{2a}p)h(p)}}=(-1)^{(p-1)/4}\sqrt p\,\ve_p^{((\f 2p)-2)(\f ap)h(p)}.
 \end{align*}
Similarly,
$$\prod_{k=1}^{n}\tan\pi\f{bk^2}p=(-1)^{(p-1)/4}\sqrt p\,\ve_p^{((\f 2p)-2)(\f bp)h(p)}.$$
If $(\f{ab}p)=-1$, then
$$\prod_{k=1}^{n}\l(\tan\pi\f{bk^2}p\times\tan\pi\f{bk^2}p\r)=\sqrt p^2=p.$$
When $(\f{ab}p)=1$, we have
$$\prod_{k=1}^{n}\l(\tan\pi\f{bk^2}p\times\tan\pi\f{bk^2}p\r)=p\ve_p^{2((\f 2p)-2)(\f ap)h(p)}.$$

Combining the above with \eqref{Tp+-}, we see that it suffices to prove \eqref{T-da} in the case $(\f{ab}p)=-1$.

Now assume $(\f{ab}p)=-1$ and set $\zeta=e^{2\pi i/p}$. By the proof of \cite [Theorem 1.1(i)]{RJ},
$T_p^{(1)}(a,b)$ is the real part of
$$D_p(a,b):=\det\l[\f{2i}{\zeta^{aj^2+bk^2}+1}\r]_{1\ls j,k\ls n},$$
and
$$D_p(a,b)=(-1)^{n/2}2^n\prod_{1\ls j<k\ls n}(\zeta^{aj^2}-\zeta^{ak^2})(\zeta^{-bj^2}-\zeta^{-bk^2}).$$
Since
$$\l(\f{a(-b)}p\r)=\l(\f{ab}p\r)=-1,$$
by Theorem \ref{Au-Th} we have
$$\prod_{1\ls j<k\ls n}(\zeta^{aj^2}-\zeta^{ak^2})(\zeta^{-bj^2}-\zeta^{-bk^2})=-\da(-ab,p)p^{(p-3)/4}
$$
and hence
$$D_p(a,b)=(-1)^{n/2}2^n\times(-1)^{n/2+1}\da(ab,p)p^{(p-3)/4}=-\da(ab,p)2^{(p-1)/2}p^{(p-3)/4}.$$
Therefore
$$T_p^{(1)}(a,b)=\Re (D_p(a,b))=-\da(ab,p)2^{(p-1)/2}p^{(p-3)/4}.$$
This proves the desired \eqref{T-da}.

By the above, we have completed our proof of Theorem \ref{Th1.4}(i). \qed

\begin{lemma} [Sun \cite{S19p}]\label{Lem-zeta} Let $p>3$ be a prime with $p\eq3\pmod4$. Let $\zeta=e^{2\pi i/p}$, and $a\in\Z$ with $p\nmid a$. Then
\begin{equation}\label{3.8}\prod_{k=1}^{(p-1)/2}(1-\zeta^{ak^2})=(-1)^{(h(-p)+1)/2}\l(\f ap\r)\sqrt {p}\,i,\end{equation}
and
\begin{equation}\label{3.9}\begin{aligned}&\ \prod_{1\ls j<k\ls(p-1)/2}(\zeta^{aj^2}-\zeta^{ak^2})
\\=&\ \begin{cases}(-p)^{(p-3)/8}&\t{if}\ p\eq3\pmod8,
\\(-1)^{(p+1)/8+(h(-p)-1)/2}(\f ap)p^{(p-3)/8}i&\t{if}\ p\eq7\pmod8,
\end{cases}\end{aligned}\end{equation}
where $h(-p)$ denotes the class number of the quadratic field $\Q(\sqrt{-p})$.
Also,
\begin{equation}\label{3.10}\prod_{1\ls j<k\ls(p-1)/2}(\zeta^{aj^2}+\zeta^{ak^2})=1,
\end{equation}
\end{lemma}

The following result can be found in \cite[Lemma 2.5]{RJ}.

\begin{lemma} [Sun \cite{RJ}]\label{Lem-S} Let $p>3$ be a prime with $p\eq3\pmod4$. Let $\zeta=e^{2\pi i/p}$, and $a,b\in\Z$ with $(\f{ab}p)=1$. Then
\begin{equation}\label{3.11}\prod_{j=1}^{(p-1)/2}\,\prod_{k=1}^{(p-1)/2}(1-\zeta^{aj^2+bk^2})
=(-1)^{(h(-p)-1)/2}\l(\f ap\r)p^{(p-1)/4}i.\end{equation}
\end{lemma}

\medskip
\noindent{\tt Proof of Theorem} 1.3(ii).  By \cite[Lemma 2.1]{RJ},
$$T_p^{(1)}(a,b,x)=c+dx$$
for some real numbers $c$ and $d$ not depending on $x$. So, it suffices to determine the value of
$T_p^{(1)}(a,b,i)$.

Let $n=(p-1)/2$ and $\zeta=e^{2\pi i/p}$. Then $\prod_{k=1}^n\zeta^{k^2}=1$ since
$$\sum_{k=0}^n k^2=\f{b(n+1)(2n+1)}6=\f{p^2-1}{24}p\eq0\pmod p.$$
For any integer $r$, clearly
$$i+\tan\pi\f rp=i+\f{(e^{i\pi r/p}-e^{-i\pi r/p})/(2i)}{(e^{i\pi r/p}+e^{-i\pi r/p})/2}
=i-i\f{\zeta^r-1}{\zeta^r+1}=\f{2i}{\zeta^r+1}.$$
Thus, with the aid of Lemma \ref{Lem-C}, we have
\begin{align*}T_p^{(1)}(a,b,i)&=\det\l[\f{2i}{\zeta^{aj^2+bk^2}+1}\r]_{1\ls j,k\ls n}
\\&=\prod_{k=1}^n\f{2i}{\zeta^{bk^2}}\times\det\l[\f1{\zeta^{aj^2}+\zeta^{-bk^2}}\r]_{1\ls j,k\ls n}
\\&=\f{2^{n}i(i^2)^{(n-1)/2}}{\zeta^{b\sum_{k=1}^{n}k^2}}
\times\f{\prod_{1\ls j<k\ls n}(\zeta^{aj^2}-\zeta^{ak^2})(\zeta^{-bj^2}-\zeta^{-bk^2})}
{\prod_{j=1}^{n}\prod_{k=1}^{n}(\zeta^{aj^2}+\zeta^{-bk^2})}
\end{align*}
and hence
\begin{equation}\label{T-mid} T_p^{(1)}(a,b,i)=i(-1)^{(p-3)/4}2^{(p-1)/2}\times\f{\prod_{1\ls j<k\ls n}(\zeta^{aj^2}-\zeta^{ak^2})(\zeta^{-bj^2}-\zeta^{-bk^2})}
{\prod_{j=1}^{n}\prod_{k=1}^{n}(\zeta^{aj^2+bk^2}+1)}.
\end{equation}

By Lemma \ref{Lem-zeta},
$$\prod_{1\ls j<k\ls n}(\zeta^{aj^2}-\zeta^{ak^2})(\zeta^{-bj^2}-\zeta^{-bk^2})
=\begin{cases} p^{(p-3)/4}&\t{if}\ p\eq3\pmod p,\\(\f{ab}p)p^{(p-3)/4}&\t{if}\ p\eq7\pmod 8.
\end{cases}$$
If $(\f{ab}p)=-1$, then $(\f{b}p)=(\f{-a}p)$ and hence
\begin{align*}\prod_{j=1}^{n}\prod_{k=1}^{n}(\zeta^{aj^2+bk^2}+1)&=\prod_{j=1}^n\prod_{k=1}^n
(\zeta^{aj^2-ak^2}+1)
=\prod_{j=1}^n\prod_{k=1}^n(\zeta^{aj^2}+\zeta^{ak^2})
\\&=\prod_{k=1}^n(2\zeta^{ak^2})\times\prod_{1\ls j<k\ls n}(\zeta^{aj^2}+\zeta^{ak^2})^2
=2^{(p-1)/2}
\end{align*}
by \eqref{3.10}.
If $(\f{ab}p)=1$, then by Lemma \ref{Lem-S} we have
\begin{align*}\prod_{j=1}^{n}\prod_{k=1}^{n}(\zeta^{aj^2+bk^2}+1)
=\prod_{j=1}^n\prod_{k=1}^n\f{1-\zeta^{2aj^2+2bk^2}}{1-\zeta^{aj^2+bk^2}}
=\f{(\f{2a}p)}{(\f ap)}=\l(\f 2p\r)=(-1)^{(p+1)/4}.
\end{align*}

Combining \eqref{T-mid} with the last paragraph, we see that if $(\f{ab}p)=-1$ then
$$c+di=T_p^{(1)}(a,b,i)=i(-1)^{(p-3)/4}2^{(p-1)/2}\times\f{(-p)^{(p-3)/4}}{2^{(p-1)/2}}=ip^{(p-3)/4}
$$
and hence
$$T_p^{(1)}(a,b,x)=c+dx=p^{(p-3)/4}x.$$
Similarly, when $(\f{ab}p)=1$ we have
$$c+di=T_p^{(1)}(a,b,i)=i(-1)^{(p-3)/4}2^{(p-1)/2}\times\f{p^{(p-3)/4}}{(-1)^{(p+1)/4}}=-i2^{(p-1)/2}p^{(p-3)/4}
$$
and hence
$$T_p^{(1)}(a,b,x)=c+dx=-2^{(p-1)/2}p^{(p-3)/4}x.$$
This concludes our proof of Theorem \ref{Th1.4}(ii). \qed

\section{Proof of Theorems \ref{Th-tan} and \ref{Th-cot}}
\setcounter{lemma}{0}
\setcounter{theorem}{0}
\setcounter{corollary}{0}
\setcounter{remark}{0}
\setcounter{equation}{0}

\medskip
\noindent {\tt Proof of Theorem \ref{Th-tan}}. Note that $\bar T_p(a,b,x)=\det[t_{jk}]_{0\ls j,k\ls n}$,
where $n=(p-1)/2$ and
$$t_{jk}=\begin{cases}1&\t{if}\ j=0,\\x+\tan\pi\f{aj^2+bk^2}p&\t{if}\ j>0.\end{cases}$$

Let $k\in\{1,\ldots,n\}$. Clearly, $t_{0k}-t_{00}=0$. Let $\zeta=e^{2\pi i/p}$. As
$$\tan\pi y=\f{2\sin \pi y}{2\cos\pi y}=\f{(e^{i\pi y}-e^{-i\pi y})/i}{e^{i\pi y}+e^{-i\pi y}}=\f{2i}{e^{2\pi i y}+1}-i$$
for all $y\in\R$ with $2y\not\in\{2m+1:\ m\in\Z\}$,
for each $j=1,\ldots,n$ we have
\begin{align*}t_{jk}-t_{j0}&=\f{2i}{\zeta^{aj^2+bk^2}+1}-\f{2i}{\zeta^{aj^2}+1}
=\f{1-\zeta^{bk^2}}{1+\zeta^{-aj^2}}
\times\f{2i}{\zeta^{aj^2+bk^2}+1}\\&=\f{(1-\zeta^{aj^2})(1-\zeta^{bk^2})}{1-\zeta^{-2aj^2}}
\times \l(i+\tan\pi\f{aj^2+bk^2}p\r).
\end{align*}

In view of the last paragraph, via all the columns (except for the first column)
of $\bar T_p(a,b,x)$ minus the first column, we see that
\begin{equation}\label{tp}\bar T_p(a,b,x)=\det[t_{jk}-t_{j0}]_{1\ls j,k\ls n}=\f{\prod_{k=1}^n(1-\zeta^{-ak^2})(1-\zeta^{bk^2})}{\prod_{j=1}^n(1-\zeta^{-2aj^2})}\times T_p^{(1)}(a,b,i).
\end{equation}

{\it Case} 1. $p\eq1\pmod 4$.

In this case, by Lemma \ref{Lem-h} we have
\begin{align*}\f{\prod_{k=1}^n(1-\zeta^{-ak^2})(1-\zeta^{bk^2})}{\prod_{j=1}^n(1-\zeta^{-2aj^2})}
&=\f{\sqrt p\,\ve_p^{-(\f {-a}p)h(p)}\sqrt p\,\ve_p^{-(\f bp)h(p)}}{\sqrt p\,\ve_p^{-(\f{-2a}p)h(p)}}
\\&=\sqrt p\,\ve_p^{((\f{2a}p)-(\f ap)-(\f bp))h(p)}=\begin{cases}\sqrt p\,\ve_p^{(\f ap)((\f 2p)-2)h(p)}\\\sqrt p\,\ve_p^{(\f{2a}p)h(p)}.\end{cases}
\end{align*}
Combining this with \eqref{tp}, \eqref{Tp+-} and Theorem \ref{Th1.4}(i), we obtain
the desired result concerning the exact value of $\bar T_p(a,b,x)$.

{\it Case} 2. $p\eq3\pmod 4$.

In this case, by Lemma \ref{Lem-zeta} we have
\begin{align*}&\ \f{\prod_{k=1}^n(1-\zeta^{-ak^2})(1-\zeta^{bk^2})}{\prod_{j=1}^n(1-\zeta^{-2aj^2})}
\\=&\ (-1)^{\f{h(-p)+1}2}\l(\f{b}p\r)\sqrt p\,i\times\f{(\f {-a}p)}{(\f {-2a}p)}=(-1)^{\f{h(-p)+1}2}\l( \f{2b}p\r)\sqrt p\,i.
\end{align*}
Combining this with \eqref{tp} and \eqref{ConjT}, we obtain the desired \eqref{Tp3}.

In view of the above, we have completed the proof of Theorem \ref{Th-tan}. \qed

\medskip
\noindent {\tt Proof of Theorem \ref{Th-cot}}.
Let $k\in\{1,\ldots,n\}$. Clearly, $c_{0k}-c_{00}=0$. Let $\zeta=e^{2\pi i/p}$. As
$$\cot\pi y=\f{2\cos\pi y}{2\sin\pi y}=\f{e^{i\pi y}+e^{-i\pi y}}{(e^{i\pi y}-e^{-i\pi y})/i}=i+\f{2i}{e^{2\pi i y}-1}\ \ \t{for all}\ y\in\R\sm\Z,$$
for each $j=1,\ldots,n$ we have
\begin{align*}c_{jk}-c_{j0}&=\f{2i}{\zeta^{aj^2+bk^2}-1}-\f{2i}{\zeta^{aj^2}-1}
=\f{1-\zeta^{bk^2}}{1-\zeta^{-aj^2}}
\times\f{2i}{\zeta^{aj^2+bk^2}-1}\\&=\f{1-\zeta^{bk^2}}{1-\zeta^{-aj^2}}
\times \l(-i+\cot\pi\f{aj^2+bk^2}p\r).
\end{align*}

In view of the last paragraph, via all the columns (except for the first column)
of $\bar C_p(a,b,x)$ minus the first column, we see that
\begin{equation}\label{cp}\bar C_p(a,b,x)=\det[c_{jk}-c_{j0}]_{1\ls j,k\ls n}=\f{\prod_{k=1}^n(1-\zeta^{bk^2})}{\prod_{j=1}^n(1-\zeta^{-aj^2})}\times C_p(a,b,-i).
\end{equation}

{\it Case} 1. $p\eq1\pmod 4$.

In this case, by Lemma \ref{Lem-h} we have
$$\f{\prod_{k=1}^n(1-\zeta^{bk^2})}{\prod_{j=1}^n(1-\zeta^{-aj^2})}
=\f{\sqrt p\,\ve_p^{-(\f bp)h(p)}}{\sqrt p\,\ve_p^{-(\f{-a}p)h(p)}}=\ve_p^{2(\f ap)h(p)}.$$
Combining this with \eqref{cp} and \eqref{cot}, we obtain
$$\bar C_p(a,b,x)=(-1)^{(p+3)/4}\da(ab,p)\f{2^{(p-1)/2}}{\sqrt p}\ve_p^{2(\f ap)h(p)}.$$

{\it Case} 2. $p\eq3\pmod 4$.

In this case, by Lemma \ref{Lem-zeta} we have
$$\f{\prod_{k=1}^n(1-\zeta^{bk^2})}{\prod_{j=1}^n(1-\zeta^{-aj^2})}
=\f{(\f bp)}{(\f {-a}p)}=\l(\f{-ab}p\r)=-1.$$
Combining this with \eqref{cp} and \eqref{C}, we obtain
$$\bar C_p(a,b,x)=(-1)^{\f{h(-p)-1}2}\l(\f ap\r)\f{2^{(p-1)/2}}{\sqrt p}.$$

In view of the above, we have completed the proof of Theorem \ref{Th-cot}. \qed

\section{Some conjectures}
\setcounter{lemma}{0}
\setcounter{theorem}{0}
\setcounter{corollary}{0}
\setcounter{remark}{0}
\setcounter{equation}{0}

Let $p$ be an odd prime, and let $d\in\Z$ with $(\f dp)=1$. For any $k=1,\ldots,(p-1)/2$, we have
$$\sum_{j=1}^{(p-1)/2}\l(\l(\f{j^2+dk^2}p\r)+\f2{p-1}\r)=-1+\f{p-1}2\times\f2{p-1}=0.$$
with the aid of \eqref{jk^2}.
Thus
$$\det\l[\l(\f{j^2+dk^2}p\r)+\f2{p-1}\r]_{1\ls j,k\ls(p-1)/2}=0,$$
and hence
\begin{equation}\det\l[x+\l(\f{j^2+dk^2}p\r)\r]_{1\ls j,k\ls(p-1)/2}=\l(1-\f{p-1}2x\r)S(d,p).
\end{equation}
by \cite[Lemma 2.1]{RJ}. Recall that $T(d,p)=\f{p-1}2S(d,p)$ by \cite[(1.20)]{S19}.
Thus, by applying \cite[Lemma 2.1]{RJ} we get that
\begin{align*}&\ \det\l[x+\l(\f{j^2+dk^2}p\r)\r]_{0\ls j,k\ls(p-1)/2}
\\=&\ T(d,p)+x\det\l[\l(\f{j^2+dk^2}p\r)-2\r]_{1\ls j,k\ls(p-1)/2}
\\=&\ \f{p-1}2S(d,p)+x\l(1-2\times\f{p-1}2\r)S(d,p).
\end{align*}
Therefore
\begin{equation}\begin{aligned}&\ \det\l[x+\l(\f{j^2+dk^2}p\r)\r]_{0\ls j,k\ls(p-1)/2}\\=&\ \l(px+\f{p-1}2\r)S(d,p)=\l(1+\f{2px}{p-1}\r)T(d,p).
\end{aligned}\end{equation}

Now let $p$ be an odd prime, and let $d\in\Z$ with $(\f dp)=-1$.
Then $S(d,p)=0$ by \cite[(1.15)]{S19}, and hence
\begin{equation}\label{xT}\det\l[x+\l(\f{j^2+dk^2}p\r)\r]_{0\ls j,k\ls(p-1)/2}=T(d,p)+S(d,p)x=T(d,p)
\end{equation}
with the aid of \cite[Lemma 2.1]{RJ}. Note that
$$1+\l(\f{0^2+d0^2}p\r)=1\ \t{and}\ \ 1+\l(\f{0^2+dk^2}p\r)=0\ \t{for all}\ k=1,\ldots,\f{p-1}2.$$
Thus
$$\det\l[1+\l(\f{j^2+dk^2}p\r)\r]_{0\ls j,k\ls(p-1)/2}=\det\l[1+\l(\f{j^2+dk^2}p\r)\r]_{1\ls j,k\ls(p-1)/2},$$
and hence
\begin{equation}\begin{aligned}&\ \det\l[x+\l(\f{j^2+dk^2}p\r)\r]_{1\ls j,k\ls(p-1)/2}
\\=&\ x\det\l[1+\l(\f{j^2+dk^2}p\r)\r]_{1\ls j,k\ls(p-1)/2}
=xT(d,p)\end{aligned}
\end{equation}
in light of \cite[Lemma 2.1]{RJ} and \eqref{xT}.

Let $p>3$ be a prime, and let $d\in\Z$ with $(\f dp)=-1$. By \eqref{d/p},
\begin{equation}\label{(p-1)/2}T(d,p)=\l(\f{p-1}2\r)^2\det\l[\l(\f{j^2+dk^2}p\r)\r]_{2\ls j,k\ls(p-1)/2}.
\end{equation}
 If $p\eq3\pmod4$, then $T(d,p)=T(-1,p)$ by \cite[(1.14)]{S19},
 and $T(-1,p)$ is an integer square by Cayley's theorem (cf. \cite[Prop. 2.2]{St})
 since it is skew-symmetric and of even order.

\begin{conjecture} \label{Conj-y} Let $p$ be a prime with $p\eq1\pmod4$.
Then, there is a positive integer $t_p$ with $(\f{t_p}p)=1$
such that for any $d\in\Z$ with $(\f dp)=-1$, we have
\begin{equation}T(d,p)=2^{(p-3)/2}\l(\f{p-1}4t_p\r)^2\,\sum_{x=1}^{(p-1)/2}\l(\f{x(x^2+d)}p\r),
\end{equation}
which has the equivalent form
\begin{equation}\det\l[\l(\f{j^2+dk^2}p\r)\r]_{2\ls j,k\ls(p-1)/2}=2^{(p-7)/2}t_p^2\,\sum_{x=1}^{(p-1)/2}\l(\f{x(x^2+d)}p\r).
\end{equation}
\end{conjecture}
\begin{remark} For any prime $p\eq1\pmod 4$ and $d\in\Z$ with $(\f dp)=-1$,
by Jacobsthal's theorem (cf. Theorem 6.2.9 of \cite[p.\,195]{BEW}) we have
$$p=\(\sum_{x=1}^{(p-1)/2}\l(\f{x(x^2+1)}p\r)\)^2+\(\sum_{x=1}^{(p-1)/2}\l(\f{x(x^2+d)}p\r)\)^2.$$
So Conjecture \ref{Conj-y} is a refinement of \cite[Conjecture 4.2(ii)]{S19}.
We have verified Conjecture \ref{Conj-y} for all primes $p<1000$ with $p\eq1\pmod4$, and found that
\begin{gather*}t_5=t_{13}=t_{17}=1,\ t_{29}=13,\ t_{37}=3^2,\ t_{41}=2\times3^2,
\\t_{53}=131,\ t_{61}=2^4\times3\times11^2,\ t_{73}=2^4\times3^3\times19\times109,
\\t_{89}=109\times199\times8273\ \t{and}\ t_{97}=2^9\times3^2\times47^2\times79.
\end{gather*}
\end{remark}

Let $p$ be an odd prime, and let $d\in\Z$ with $(\f {-d}p)=-1$. For the matrix $A_p=[a_{jk}]_{0\ls j,k\ls (p-1)/2}$ with
$$a_{jk}=\begin{cases}1&\t{if}\ j=0,\\1/(j^2+dk^2)&\t{if}\ j>0,
\end{cases}$$ we have
\begin{align*}
\det A_p&=(-d)^{(p-1)/2}\det\l[\f1{j^2+dk^2}\r]_{1\ls j,k\ls(p-1)/2}
\\&\eq-\det\l[\f1{j^2+dk^2}\r]_{1\ls j,k\ls(p-1)/2}\pmod p;
\end{align*}
this can be seen by considering each column (except the first column) minus the first column and noting that
$$\f1{j^2+dk^2}-\f1{j^2+d0^2}=\f{-dk^2}{j^2(j^2+dk^2)}\ \ \t{for all}\ j,k=1,\ldots,\f{p-1}2.$$
Thus, with the aid of \eqref{dk^2}, we get
\begin{equation}\det A_p\eq\begin{cases}-d^{(p-1)/4}\pmod p&\t{if}\ p\eq1\pmod4,
\\(-1)^{(p-3)/4}\pmod p&\t{if}\ p\eq3\pmod4,\end{cases}
\end{equation}
and hence
\begin{equation}\l(\f{\det A_p}p\r)=(-1)^{\lfloor(p-3)/4\rfloor}=\l(\f{-2}p\r).
\end{equation}

\begin{conjecture}\label{Conj-p-2} Let $p$ be a prime with $p\eq1\pmod4$, and let $d\in\Z$ with $(\f dp)=-1$.
Then
\begin{equation}\label{p-2*}3\bar S_{p-2}(1,p)\eq S_{p-2}(1,p)\eq2\da(d,p)\sum_{x=1}^{(p-1)/2}\l(\f{x(x^2+d)}p\r)\pmod p,
\end{equation}
where $\bar S_{p-2}(1,p)=\det[s_{jk}]_{0\ls j,k\ls (p-1)/2}$ with
$$s_{jk}=\begin{cases}1&\t{if}\ j=0,\\(j^2+k^2)^{p-2}&\t{if}\ j>0.\end{cases}$$
\end{conjecture}
\begin{remark} Let $p\eq1\pmod4$ be a prime, and write $p=x^2+y^2$ with $x,y\in\Z^+$ and $2\mid y$.
Then, for any $d\in\Z$ with $(\f dp)=-1$, we have $\sum_{x=1}^{(p-1)/2}\l(\f{x(x^2+d)}p\r)=\pm y$
by Jacobsthal's theorem.
Let $q=\f{p-1}2!$. Then $(y/x)^2\eq-1\eq q^2\pmod p$ and hence
$$\l(\f yp\r)=\l(\f{qx}p\r)=\l(\f qp\r)\l(\f px\r)=\l(\f 2p\r)$$
with the aid of \cite[Lemma 2.3]{S19}.
Thus Conjecture \ref{Conj-p-2} implies that
\begin{equation}\l(\f{S_{p-2}(1,p)}p\r)=\l(\f{3\bar S_{p-2}(1,p)}p\r)=1.
\end{equation}
\end{remark}

Let $m,n\in\Z^+$ with $n$ odd. For the determinant
\begin{equation}D^{(m)}_n:=\det\l[(j^2-k^2)^m\l(\f{j^2-k^2}n\r)\r]_{1\ls j,k\ls(n-1)/2},
\end{equation}
clearly
\begin{align*}D^{(m)}_n&=\det\l[(k^2-j^2)^m\l(\f{k^2-j^2}n\r)\r]_{1\ls j,k\ls(n-1)/2}
\\&=\l((-1)^m\l(\f{-1}n\r)\r)^{(n-1)/2}D^{(m)}_n=(-1)^{(m-1)(n-1)/2}D^{(m)}_n,
\end{align*}
and hence $D^{(m)}_n=0$ when $2\mid m$ and $4\mid n-3$.
If $2\nmid m$ and $4\mid n-1$, then $D^{(m)}_n$ is skew-symmetric
and of even order, hence it is an integer square by Cayley's theorem.

\begin{conjecture} For any prime $p\eq1\pmod4$, we have
\begin{equation}\label{Dp1}\l(\f{\sqrt{D_p^{(1)}}}p\r)=(-1)^{|\{0<k<\f p4:\ (\f kp)=-1\}|}\l(\f p3\r).
\end{equation}
\end{conjecture}
\begin{remark} We have verified \eqref{Dp1} for all primes $p<1000$ with $p\eq1\pmod4$.
\end{remark}

\begin{conjecture} For any prime $p\eq1\pmod4$, we have
\begin{equation}\label{Dp3}\l(\f{\sqrt{D_p^{(3)}}}p\r)=(-1)^{|\{0<k<\f p4:\ (\f kp)=-1\}|}\l(\f p{4+(-1)^{(p-1)/4}}\r).
\end{equation}
\end{conjecture}
\begin{remark} We have verified \eqref{Dp3} for all primes $p<1000$ with $p\eq1\pmod4$.
\end{remark}

\begin{conjecture} For any positive odd integer $m$, the set
$$E(m)=\l\{p:\ p\ \t{is a prime with}\ 4\mid p-1\ \t{and}\ p\mid D_p^{(m)}\r\}$$
is finite. In particular,
\begin{gather*}E(5)=\{29\},\ E(7)=\{13,\,53\},\ E(9)=\{13,\,17,\,29\},
\ E(11)=\{17,\,29\}.
\end{gather*}
\end{conjecture}
\begin{remark} This is based on our computation. For $m=5,7,9,11$, we find those primes $p<1000$ in $E(m)$ via {\tt Mathematica}. It seems that $\{p\in E_{13}:\ p<1000\}==\{17,\,109,\,401\}.$
\end{remark}

\setcounter{conjecture}{0} \end{document}